\documentclass[11pt,reqno]{amsart}

\usepackage{geometry}
\usepackage{amsmath, amsthm, amssymb, bm, color, framed, graphicx, hyperref, mathrsfs}
\usepackage{cite}

\usepackage[all]{xy}
\usepackage{bbm}
\usepackage{cancel}
\usepackage{esint}
\definecolor{shadecolor}{RGB}{241, 241, 255}
\newtheorem{theorem}{Theorem}

\def \M {\mathcal{M}}

\def\al{\aligned}
\def\eal{\endaligned}
\def\be{\begin{equation}}
\def\ee{\end{equation}}
\def\befr{\begin{frame}}
\def\enfr{\end{frame}}
\def\lab{\label}

\def\e{\epsilon}
\def\n{\nabla}

\def\R{{\bf R}}
\def\M{{\bf M}}

\def\al{\aligned}

\def\De{\Delta}

\def\pd{\partial}

\begin{document}
\title{An improved Hamilton matrix estimates for the heat equation}
\author{Lang Qin and Qi S. Zhang}
\address{School of Mathematics, Fudan University, Shanghai, China and Department of
Mathematics  University of California, Riverside, CA 92521, USA }
\email{22110180037@m.fudan.edu.cn and qizhang@math.ucr.edu}
\date{}
\maketitle
\newcommand\dif{\mathop{}\!\mathrm{d}}
\newcommand\diver{\mathop{}\!\mathrm{div}}
\newcommand\vol{\mathop{}\!\mathrm{vol}}
\newcommand\Ric{\mathop{}\!\mathrm{Ric}}
\begin{abstract}
In this paper, we remove the assumption on the gradient of the Ricci curvature in Hamilton's matrix Harnack estimate for the heat equation on all closed manifolds, answering a question which has been around since the 1990s. New ingredients include a recent sharp Li-Yau estimate,  construction of a  suitable vector field and various use of integral arguments, iteration  and a little tensor algebra.
\end{abstract}

\section{Introduction}

 Let $\M$ be a $n$ dimensional,  compact Riemmannian manifold without boundary and $u: \M\times [0,\infty)\to\R$ be a positive solution to the heat equation
\begin{equation}
	\frac{\partial u}{\partial t}=\Delta u,
	\label{01}
\end{equation}

In 1993 R. Hamilton proved the following matrix Harnack inequality

\begin{theorem} (\cite{Ha93})
\lab{thham93}
 For any positive solution $u$ to \eqref{01}, there exist constants $B$ and $C$ depending only on the geometry of $\M$ (in particular the diameter, the volume, and the curvature and covariant derivative of the Ricci curvature) such that if $t^{n/2}u\leq B$ then
\begin{equation}\label{eq Hamilton matrix general case}
\n_i \n_j \ln u  +\frac{1}{2t}g_{ij} +C \left(1+\ln\left(\frac{B}{t^{n/2}u}\right)\right)g_{ij} \geq 0.
\end{equation}

Moreover $C=0$ if $\M$ has nonnegative sectional curvature and parallel Ricci curvature and the inequality holds for all positive solutions $u$.

\end{theorem}

Here $\n_i \n_j$ is the Heissian operator.
This result is an extension, to the matrix form,  of the Li-Yau estimate
in 1986 \cite{LY:1}, when Peter Li and S.-T. Yau developed a log gradient estimate for positive solutions $u=u(x, t)$ of the heat equation on in \cite{LY:1}.  When the Ricci curvature is non-negative, the Li-Yau  estimates  takes the form\begin{equation}
	\Delta\ln u+\frac{n}{2t}\ge 0,
	\label{03}
\end{equation}
with the equality achieved on the heat kernel on the Euclidean space $\R^n$.

  As mentioned, under the stronger assumptions that $\M$ is Ricci parallel and has weakly positive sectional curvatures.
Hamilton's estimate  becomes
\begin{equation}
	\nabla_i\nabla_j\ln u+\frac{1}{2t} g_{ij}\ge 0.
	\label{04}
\end{equation}
 Notice that the trace of (\ref{04}) is (\ref{03}).
 The estimates by Li-Yau and Hamilton and their generalizations in various settings have become powerful tools for studying the geometry, topology and analytical properties of manifolds, see the book \cite{Lib}  and survey papers, \cite{B.Chow},  \cite{Z24} e.g.

Later, in \cite{CH}, Chow and Hamilton extended (\ref{04}) to constrained case under the same assumptions and discovered new linear Harnack estimates. In \cite{CN},   Cao and Ni proved the matrix Li-Yau-Hamilton estimates on complete K\"ahler manifolds with nonnegative bisectional curvature\begin{equation}
	\nabla_\alpha\nabla_{\bar\beta}\ln u+\frac{1}{2t} g_{\alpha\bar\beta}\ge 0.
	\label{05}
\end{equation}
They removed the assumption
that the Ricci tensor being parallel in this case.

For more results about the generalization of matrix Li-Yau-Hamilton estimates in different settings, we refer the  reader to \cite{P22}\cite{R14}\cite{R19}\cite{Ch14} . In \cite{LZ} Li and Zhang proved a matrix Li-Yau-Hamilton estimate for the heat equation under the Ricci flow. As applications, they establish the monotonicity  of parabolic frequencies up to some correction factors. They also removed the extra log term $\ln\left(\frac{B}{t^{n/2}u}\right)$ in Hamilton's estimate \eqref{eq Hamilton matrix general case}.  Similar estimates are established under K\"ahler Ricci flow in
\cite{Li23}. See also  a sharp estimate in \cite{Ni05} by Ni for a different equation.

However, in general, the assumption on $\nabla Ric$ in Theorem \ref{thham93} has been considered somewhat less ideal over the years, c.f. \cite{CN} e.g. The  purpose of the paper is to remove this assumption.

Our main result is the following matrix Hamilton estimates:
\begin{theorem}
\lab{thmain}
	Let $(\M, g)$ be a closed  n-manifold, $n\ge 2$, $u:\M \times (0,T)\to\R$ be a positive solution to the heat equation \eqref{01}. Suppose there exist positive constants $K, v_0>0$ such that sectional curvatures satisfy $|sec| \le  K$, the volume of unit balls satisfy $\vol(B(p,1))\ge v_0, \forall p\in \M$.
	Then\begin{equation}
		\nabla_i\nabla_j\ln u+ \frac{\gamma(T,n,K,v_0)}{t}g_{ij}\ge 0.
		\label{06}
	\end{equation}
	\label{main}
Here $\gamma(T,n,K,v_0)$ is a non-negative constant depending only on the dimension of the manifold, the curvature bound $K$ and the lower volume bound of the unit ball $v_0$ and  the diameter of $M$.
\end{theorem}


Note that within the proof of Corollary 1.2 of \cite{Ha93} p117, Hamilton used  the assumption that the unit balls have a definite lower volume bound. This corollary was used in the proof  Corollary 4.2, and then Theorem 4.3 there and his Theorem \ref{thham93} here. Therefore the assumption of non-collapsedness in our Theorem \ref{main} was also used there. In this sense Theorem \ref{main} is a proper improvement of Theorem \ref{thham93}.

 In the general noncompact case, one cannot prove a bound like \eqref{06}. After taking trace, it is impossible to have a finite $\gamma(T,n,K,v_0)$   for a fixed $t$ even for the hyperbolic space, see \cite{Z21} for detail.

Let us explain the idea of the proof for Theorem \ref{thmain}. As in \cite{Ha93}, the starting point is the equation \eqref{14} for the tensor $-t \nabla_i\nabla_j\ln u$.
 Instead of using the tensor maximum principle, in order to remove the condition on the gradient of the Ricci curvature, a natural approach is to use integral method. So the first attempt would be to consider the equation for the norm of the Hessian of $-t \ln u$. i.e. $t |-\nabla_i\nabla_j\ln u|$. This is the usual semigroup domination method in dealing with system of equations.  However, this method is doomed to fail here since if it could prove Theorem \ref{thmain}, it would also prove that $t |\Delta \ln u| \le C$. The latter  is false even on the 2 sphere $S^2$ since it would imply, by choosing $u$ to be the heat kernel on $S^2$ and using standard small time asymptotes, that $ |\Delta r^2| \le 4C$ away from the cut-locus. Here $r$ is the distance function from the north pole. Then it would infer that $|\Delta r| \le (2C-1)/r$. However, by direct computation, we know $\Delta r= \cos r/\sin r$ for $r \in (0, \pi)$. We have reached a contradiction when $r \to \pi$. This example shows that  $-t \nabla_i\nabla_j\ln u$ may not be bounded from below.

 Instead our method is to apply the tensor equation \eqref{14} on a certain vector field generated by the eigenvector of the most positive eigenvalue of $-t \nabla_i\nabla_j\ln u$, obtaining a nonlinear scalar equation with many cross terms. Due to the complexity of the curvature terms, certain careful calculation involving linear algebra is needed to avoid using the most negative eigenvalue of $-t \nabla_i\nabla_j\ln u$, which may be unbounded as mentioned above.
Other new ingredients in the proof are the use of integral arguments  and iteration, together with the recent sharp Li-Yau bound on compact manifolds when the Ricci curvature changes sign. The nonlinear term, having a good sign, plays a crucial role.

We conclude the introduction with a list of common notations and symbols in the paper.
We will use $u=u(x, t)$ to denote a solution to the heat equation $\Delta u - \partial_t u=0$ where $\Delta$ is the Laplace-Beltrami operator; $\nabla f$ denotes the gradient  and $\nabla^2 f$  or $\nabla_i \nabla_j  f$ denote the Hessian of a scalar function $f$ respectively. Given a tensor field $T$ and vector $e_i$, $\nabla_{e_i} T$ denotes the covariant derivative in the direction of $e_i$ which may also be denoted by $T_{; i}$ in component form. $\M$ denotes a closed $n$ dimensional Riemannian manifold; $sec$ and $Ric$ denote the sectional and Ricci curvature respectively. Given a tensor field $T=T(x, t)$, we use $|T(x, t)|$ or $|T|$ to denote it is norm at $(x, t)$. Constants are denoted by $C$ with or without index. Generic constants may change in value from line to line.

\section{proof of Theorem \ref{main}}
\begin{proof}

We divide the proof into three steps.

{\it Step 1. Set up a nonlinear partial differential inequality.}

First, we choose $u=u(x, t)$ to be the heat kernel $G(x,t+\delta;y)$ for a fixed $y \in \M$ and a small $\delta>0$ and will eventually let $\delta \to 0$.
For a fixed time $T$, since $\M\times [0,T]$ is compact, the maximum of the largest eigenvalue of $-t\nabla^2\ln u$ on $\M\times [0,T]$ is reached at some point $(x_0,t_0)$.

We denote the largest eigenvalue of $-t\nabla^2\ln u $ on $\M\times [0,T]$ as $Q_0$. If $Q_0\le 0$, then $\nabla^2\ln u \ge 0$ on $\M\times [0,T]$. If this holds for all $T>0$, the proof is done. So we assume for some fixed time $T$,  we have $Q_0> 0$.
 Let $\xi_0$ be a unit eigenvector of $-t\nabla^2\ln u$  with respect to $Q_0$ at $(x_0,t_0)$.

From the results of Cheng-Li-Yau in \cite{CLY} and Cheeger-Gromov-Taylor in \cite{CGT}, we know the injectivity radius of $\M$ is bounded below by some constant  depending only on $n,K$ and $v_0$. We will denote it  as $i_0$.

We next construct a harmonic coordinate system on the ball $B(x_0,i_0)$ using Anderson's result on the lower bound of harmonic radius.

Let  $g(y)$ represent the Riemannian metric tensor at the point $y\in \M$  and fix some constants $C_1>1$, $0<\alpha<1$,
applying the no boundary version of {\it Main Lemma 2.2} in \cite{An},
 there is a small $\varepsilon_0=\varepsilon_0(K,n,C_1,\alpha)$ with the following property:
For any
 $\varepsilon(x_0)\ge \varepsilon_0 i_0 $, there is a harmonic coordinate system $\{x^i\}_{i=1}^n$ defined on $B(x_0,\varepsilon(x_0))$ such that
\begin{equation}
\al
&g(x_0)(\nabla x^i,\nabla x^j)=\delta^{ij},\\
&C_1^{-1} I\le g(y)\le C_1 I,\\
&\varepsilon(x_0)^{1+\alpha}\|g(\cdot )(\nabla x^i,\nabla x^j)\|_{C^{1,\alpha}(B(x_0,  \varepsilon(x_0) \not i_0))}\le C_1,
\eal\label{07}
\end{equation}
for all $y\in B(x_0,\varepsilon(x_0))$.   Since $\varepsilon(x_0)$ can be written as  a constant multiple of $i_0$, for convenience, we will still use $i_0$ to denote $\varepsilon(x_0)$ in the rest of the proof.

At the point $(x_0,t_0)$, the unit eigenvector $\xi_0$ defined earlier can be written as a unique linear composition of $\{\nabla x^i\}$, say $\xi_0=a_i\nabla x^i(x_0)$, $a_i\in\mathbb\R$. From (\ref{07}) and $|\xi_0|=1$, we know that
\begin{equation}
	\sum\limits_{i=1}^n a_i^2=1.
	\label{08}
\end{equation}

Define \begin{equation}
	F(y):=a_ix^i(y),\qquad \beta(y):=\nabla F(y), \qquad \forall y\in B(x_0,i_0).
	\label{09}
\end{equation}

We use the Cheeger-Colding cutoff function, see \cite{CC}
\begin{equation*}
\eta=\left\{\al
&1, &\mbox{on}\quad &B\left(x_0,\frac{i_0}{2}\right)\\
&0, &\mbox{on}\quad &\M\setminus B\left(x_0, \frac{3i_0}{4}\right)
\eal\right.  ,
\end{equation*} which is a non-negative smooth function on $\M$ satisfying \begin{equation}
0\le\eta\le 1,\quad |\nabla \eta|\le \frac{C_2}{i_0}\eta^{1-\delta}, \quad |\Delta\eta| \le\frac{C_2}{i_0^2}.\label{10}\end{equation}

Define  the vector field
\begin{equation}
\al
	\xi:=\begin{cases}
\eta\beta, \quad  \mbox{on}\quad B(x_0,i_0)\\
0, \quad \mbox{on}\quad \M\setminus B(x_0,i_0).
\end{cases}
\eal
\label{11}
\end{equation}
 Using (\ref{07}) and (\ref{08}) ,we calculate, $\forall y\in B(x_0,i_0)$,
 that \begin{equation}
 	\sqrt{C_1^{-1}}= \left(C_1^{-1} \sum\limits_{i=1}^n a_i^2\right)^{\frac{1}{2}}\le  |\beta|=g(y)( a_i\nabla x^i, a_j\nabla x^j)^{\frac{1}{2}}\le \left(C_1 \sum\limits_{i=1}^n a_i^2\right)^{\frac{1}{2}} =\sqrt{C_1}, 	 \label{12}
 \end{equation}
Therefore, $\xi$ is a well-defined smooth vector field on $\M$, $\xi(x_0)=\xi_0$ and $|\xi|\le \sqrt{C_1}$.

We will bound the  norms $\||\nabla\xi|\|_\infty$ and $\||\Delta \xi|\|_\infty$ on $B(x_0,i_0)$ using (\ref{07}) and (\ref{10}).  Observe that
\begin{equation}\al
	&\nabla_{\frac{\pd}{\pd x^k}} \xi=\frac{\pd\eta}{\pd x^k}\beta+\eta\nabla_{\frac{\pd}{\pd x^k}}\beta,\\
	&|\nabla\xi|=\left(\left\langle \nabla_{\frac{\pd}{\pd x^k}} \xi, \nabla_{\frac{\pd}{\pd x^l}} \xi\right\rangle g^{kl}\right)^{\frac{1}{2}},\\
	& \| |\nabla\xi| \|_\infty\le C_1 C_2i_0^{-1}+\sqrt{C_1}\|\nabla\beta\|_\infty.  \\
	&\De\beta=\De \nabla F=\nabla \De F+Ric(\nabla F)=  Ric (\beta),\\
	&|\De \xi| \le  |\De\eta| |\beta|+|\eta||  Ric(\beta)|+2|\nabla\eta||\nabla\beta|,\\
	&\| |\De\xi| \|_\infty\le C_2i_0^{-2}\sqrt{C_1}+(n-1)K\sqrt{C_1}+2C_2i_0^{-1}\|\nabla\beta\|_\infty.
	\eal\label{13}
	\end{equation}Here we used the Bochner's formula and $\Delta F=0$ from (\ref{10}) and $ Ric$ denotes the Ricci curvature. The equation $\Delta \beta = Ric(\beta)$ stands for
\[
<\Delta \beta, X> = Ric(\beta, X)
\]for any smooth vector field $X$.
	 We can apply Theorem 3.2 from \cite{WZ} to prove $\|\nabla \beta\|_\infty$ is finite. Note that theorem was stated for eigen-forms of the Hodge Laplacian. But the proof works verbatim for the equation $\Delta \beta =Ric(\beta)$. The point is that one only needs the condition that the sectional curvature being bounded not the gradient of the Ricci curvature being bounded.

 We use $H$ to represent the $(0,2)$-tensor $-t\nabla^2\ln u$.
At any point  $y\in M$, we can choose an orthonormal coordinates $\{e_i\}$, and denote $H_{ij}=H(e_i,e_j)$. According to the basic calculation, see \cite{Ha93},
\begin{equation}
	\begin{split}
		(\Delta-\partial_t)H_{ij}=\frac{2}{t}H_{ij}^2-\frac{1}{t}H_{ij}
		-2R_{ikjl}H_{kl}+R_{ik}H_{kj}+R_{jk}H_{ki}\\+2tR_{ikjl}(\nabla^k\ln u)(\nabla^l\ln u)-2(\nabla_kH_{ij})(\nabla^k\ln u)\\+t(\nabla_l R_{ij}-\nabla_i R_{jl}-\nabla_j R_{il})(\nabla^l \ln u).
	\end{split}
	\label{14}
\end{equation} Here $H^2_{ij} = H_{ik} H_{kj}$.

The equation \eqref{14} is for tensor. After plugging two vectors $X,Y$ into it, we have
\begin{equation}\al
&	(\De H-\partial_tH)(X,Y)\\
&=\frac{2}{t}\langle H(X,\bullet),H(Y,\bullet)\rangle -\frac{1}{t}H(X,Y)-2\langle Rm(X,\bullet,Y,\bullet),H\rangle\\
&+\langle  Ric(X,\bullet), H(Y,\bullet)\rangle
	+\langle  Ric(Y,\bullet), H(X,\bullet)\rangle \\
&+2tRm(X,\nabla\ln u,Y,\nabla\ln u)-2\langle(\nabla H)(\bullet,X,Y),\nabla\ln u\rangle\\
&	+t\langle (\nabla   Ric)(\bullet,X,Y),\nabla\ln u\rangle - t\langle (\nabla   Ric)(X,\bullet,Y),\nabla\ln u\rangle -t\langle (\nabla Ric)(Y,\bullet,X),\nabla\ln u\rangle.
	\eal\label{15}
\end{equation}
Notice the $\langle,\rangle$ here means the inner product of tensor induced by the metric $g$. For instance, if $\{ e_i \}$ is an orthonormal basis for $T_x \M$ at a point $x$, then $\langle H(X,\bullet),H(Y,\bullet)\rangle = \sum^n_{i=1} H(X,e_i) H(Y, e_i)$, and $\langle (\nabla  Ric)(\bullet,X,Y),\nabla\ln u\rangle= \sum^n_{i=1} (\nabla_{e_i}   Ric)(X,Y)\nabla_{e_i} \ln u$ etc.

Now we define the scalar function
\be
\lab{defQ}
Q:=H(\xi,\xi)
\ee on $\M$.
By the definition of $\xi$, the function $Q$ vanishes outside $B(x_0,i_0)$. In $B(x_0,i_0)$, we use the harmonic coordinates $\{x_i\}_{i=1}^n$ to be the local chart, $H=Q_{ij}\dif x^i\otimes \dif x^j$, $\xi=\xi^i\frac{\pd}{\pd x^i}$,$ \nabla H=Q_{ij;k}\dif x^k\otimes \dif x^i\otimes \dif x^j$, $\nabla_k \xi=\xi^i_{;k}\frac{\pd}{\pd x^i}$, $Q=\xi^iQ_{ij}\xi^j$.

Using \eqref{15},  we establish the equation for $Q$:
\begin{equation}\al
  (\Delta-\partial_t)Q	&=\De (H(\xi,\xi))-\pd_t(H(\xi,\xi))
  \\&=(\De H-\pd_t H)(\xi,\xi)+2H(\De\xi,\xi)+4(\nabla H)(\nabla \xi,\xi)+2H(\nabla \xi,\nabla \xi)
  \\&=\frac{2}{t}\langle H(\xi,\bullet),H(\xi,\bullet)\rangle -\frac{1}{t}H(\xi,\xi)-2\langle Rm(\xi,\bullet,\xi,\bullet),H\rangle+\langle  Ric(\xi,\bullet), H(\xi,\bullet)\rangle\\&
	\quad+\langle  Ric(\xi,\bullet), H(\xi,\bullet)\rangle +2tRm(\xi,\nabla\ln u,\xi,\nabla\ln u)-2\langle(\nabla H)(\bullet,\xi,\xi),\nabla\ln u\rangle\\&
	\quad+t\langle (\nabla   Ric)(\bullet,\xi,\xi),\nabla\ln u\rangle - t\langle (\nabla  Ric)(\xi,\bullet,\xi),\nabla\ln u\rangle -t\langle (\nabla   Ric)(\xi,\bullet,\xi),\nabla\ln u\rangle\\
	&\quad+2H(\De\xi,\xi)+4(\nabla H)(\nabla \xi,\xi)+2H(\nabla \xi,\nabla \xi)
  \\  &= \frac{2}{t} \xi^iQ_{ij}^2\xi^j-\frac{1}{t}Q
		-2R_{ikjl}Q_{pq}g^{pk}g^{lq}\xi^i\xi^j+2\xi^i R_{ik}Q_{kj}\xi^j
\\&\quad+2tR_{ikjl}\xi^i\xi^j(\nabla^k\ln u)(\nabla^l\ln u) -2(\nabla_k Q-2Q_{ij}\xi^i \xi^j_{;k})(\nabla^k\ln u)
\\&\quad+t( R_{ij;l}- R_{jl;i}-R_{il;j})(\nabla^l \ln u)\xi^i\xi^j +2Q_{ij}(\Delta \xi)^i\xi^j+4 Q_{ij;l} g^{kl}\xi^i_{;k}\xi^j+2Q_{ij}\xi^i_{;k}\xi^j_{;l}g^{kl}.
	 \eal\label{16}
	 \end{equation}
Here the notations $(\nabla H)(\nabla \xi ,\xi)=\sum\limits_{i=1}^{n}(\nabla_{e_i} H)(\nabla_{e_i} \xi ,\xi),\quad H(\nabla \xi,\nabla \xi)=\sum\limits_{i=1}^{n}H(\nabla_{e_i}\xi, \nabla_{e_i} \xi)$ are used.

	 We now handle the terms with $R_{ikjl}$ on the right hand side of (\ref{16}).
Using $|sec|\le K$, we have \begin{equation}
	 \begin{aligned}
	 	&2tR_{ikjl}\xi^i\xi^j(\nabla^k\ln u)(\nabla^l\ln u)\\
	 	=&2t\mathrm{Rm}(\xi,\nabla\ln u,\xi,\nabla\ln u)\\
	 	\ge& -2tK|\xi|^2|\nabla\ln u|^2.
	 \end{aligned}
 	\label{17}
 \end{equation}

	Next,
		 \begin{equation}
	 	\begin{aligned}
	 	&-2R_{ikjl}Q_{pq}g^{pk}g^{lq}\xi^i\xi^j\\
	 		 =&-2(R_{ikjl}+K(g_{ij}g_{kl}-g_{il}g_{jk}))Q_{pq}g^{pk}g^{lq}\xi^i\xi^j+2Kg^{pq}Q_{pq}|\xi|^2-2KQ\\
	 		 \ge&-2(R_{ikjl}+K(g_{ij}g_{kl}-g_{il}g_{jk}))Q_0 g^{kl}\xi^i\xi^j-2K|\xi|^2t\Delta\ln u-2KQ\\
	 		 \ge& -4KQ_0(g_{ij}g_{kl}-g_{il}g_{jk})g^{kl}\xi^i\xi^j-2K|\xi|^2t\Delta\ln u-2KQ \\
	 		 =&-4KQ_0(n-1)|\xi|^2-2K|\xi|^2 t\Delta\ln u-2KQ\\
	 		 \ge &-2K|\xi|^2 t\Delta\ln u-4(n-1)C_1KQ_0-2C_1KQ_0,
	 	\end{aligned}\label{18}
	 \end{equation}
where we used the fact that $g^{kp} Q_{pq} g^{ql}\le Q_0 g^{kp} g_{pq} g^{ql}=Q_0g^{kl}$ by the definition of $Q_0$ and positive definiteness of $g^{ij}$ and the fact that $[R_{ikjl}+K(g_{ij}g_{kl}-g_{il}g_{jk})]\xi^i\xi^j\ge 0$ by $|sec|\le K$ in the first inequality,	  $|sec|\le K$ in the second one, and $|\xi|^2\le C_1$, $Q\le Q_0|\xi|^2$ in the last one.

Substituting (\ref{17}) and (\ref{18}) into (\ref{16}), we have
\begin{equation}
\al	
&	 (\Delta-\partial_t)Q+2(\nabla_k Q-2Q_{ij}\xi^i \xi^j_{;k})(\nabla^k\ln u)\\
&\ge \frac{2}{t}\xi^iQ^2_{ij}\xi^j-\frac{1}{t}Q-2K|\xi|^2 t|\nabla\ln u|^2-2K|\xi|^2 t\Delta\ln u\\
&\qquad -2(2n-1)C_1KQ_0+2\xi^i R_{ik}Q_{kj}\xi^j+t( R_{ij;l}- R_{jl;i}-R_{il;j})(\nabla^l \ln u)\xi^i\xi^j\\
 &\qquad +2Q_{ij}(\Delta \xi)^i\xi^j+4 Q_{ij;l} g^{kl}\xi^i_{;k}\xi^j+2Q_{ij}\xi^i_{;k}\xi^j_{;l}g^{kl}.
\label{19}
\eal
\end{equation}
	
Using the equality $- \Delta \ln u = |\nabla \ln u |^2 - \partial_t \ln u$, we can cancel the term involving $|\nabla  \ln u|^2$ on the right hand side.
 By definition, $Q_{ij}\equiv 0$ outside $B(x_0,i_0)$,
 so $Q_+=\max\{Q,0\}$ is a subsolution of the following inequality in the weak sense on $\M \times (0,t_0)$,\begin{equation}\begin{split}	
(\Delta-\partial_t)Q_++2(\nabla_k Q-2Q_{ij}\xi^i \xi^j_{;k})(\nabla^k\ln u)\ge \frac{2}{t}\xi^iQ^2_{ij}\xi^j-\frac{1}{t}Q_+-2|\xi|^2K t\partial_t \ln u
\\-2(2n-1)C_1KQ_0+2\xi^i R_{ik}Q_{kj}\xi^j+t( R_{ij;l}- R_{jl;i}-R_{il;j})(\nabla^l \ln u)\xi^i\xi^j\\ +2Q_{ij}(\Delta \xi)^i\xi^j+4 Q_{ij;l} g^{kl}\xi^i_{;k}\xi^j+2Q_{ij}\xi^i_{;k}\xi^j_{;l}g^{kl}
 .
	 \end{split}\label{20}
	 \end{equation}
In the next two steps, we will use an integral argument inspired by Moser's iteration to prove that $Q$ is bounded from above. Due to the complexity of the curvature terms, additional work is needed. The main thing is to avoid the most negative eigenvalue of $Q_{ij}$,  which, as mentioned, may be unbounded.
\medskip

{\it Step 2. $L^m$ bound for a fixed large number $m$.}

We begin with two derivative estimates of the heat kernel $u=G(x, t+\delta; y)$ essentially due to Hamilton, which will be needed in the sequel.
	
 Using Hamilton's gradient bound for  $|\nabla \ln u|$ (Theorem 1.1) and  $\partial_t \ln u$ (Lemma 4.1) in \cite{Ha93}, we can obtain the following estimates
\begin{equation}
		t|\nabla\ln u|^2\le 2(1+(n-1)Kt)\left(C(n)(1+K+Kt)+\frac{D_0^2}{2t}\right),
\label{21}
\end{equation}
\begin{equation}
		t \partial_t \ln u\le 2(1+(n-1)Kt)\left(C(n)(1+K+Kt)+\frac{D_0^2}{2t}\right);
\label{22}
\end{equation}
where $D_0$ is the diameter of $\M$, $C(n)$ is some dimensional constant.
The proof will be outlined in the appendix for completeness.

We plug  (\ref{22}) into the third term on the right hand side of (\ref{20}) and   denote, for simplicity,  $C_3=2(1+(n-1)KT)(1+K+KT)C(n)$ for some $T>t_0$ fixed. Then we deduce
\begin{equation}	
\al
&(\Delta-\partial_t)Q_++2(\nabla_k Q-2Q_{ij}\xi^i \xi^j_{;k})(\nabla^k\ln u)\\
&\ge \frac{2}{t}\xi^iQ^2_{ij}\xi^j-\frac{1}{t}Q_+-\frac{2C_1C_3D_0^2K}{t}
	 -2C_1C_3K
\\
&-2(2n-1)C_1KQ_0+2\xi^i R_{ik}Q_{kj}\xi^j+t( R_{ij;l}- R_{jl;i}-R_{il;j})(\nabla^l \ln u)\xi^i\xi^j\\
&+2Q_{ij}(\Delta \xi)^i\xi^j+4 Q_{ij;l} g^{kl}\xi^i_{;k}\xi^j+2Q_{ij}\xi^i_{;k}\xi^j_{;l}g^{kl}.
	 \label{23}
\eal
	 \end{equation}

We choose $ t Q_+^{2m}$ as the test function for (\ref{23}).
Then we arrive at
\begin{equation}\begin{aligned}
	&\int_{0}^{t_0}\int_{\M} 2\xi^iQ^2_{ij}\xi^jQ_+^{2m}\dif x\dif t\\
\le& \underbrace{\int_{0}^{t_0}\int_{\M} t Q_+^{2m} (\Delta-\partial_t)Q_+ \dif x\dif t}_I\\
&+\underbrace{2\int_{0}^{t_0}\int_{\M} 2(\nabla_k Q-2Q_{ij}\xi^i \xi^j_{;k})(\nabla^k\ln u) tQ_+^{2m}\dif x\dif t}_{II}\\
&+\underbrace{\int_{0}^{t_0}\int_{\M} Q_+^{2m+1}+2C_1C_3D_0^2K Q_+^{2m} +2C_1C_3Kt Q_+^{2m} +2(2n-1)C_1K tQ_0 Q_+^{2m} \dif x\dif t}_{III} \\
&\underbrace{-\int_{0}^{t_0}\int_{\M} [2\xi^i R_{ik}Q_{kj}\xi^j+t( R_{ij;l}- R_{jl;i}-R_{il;j})(\nabla^l \ln u)\xi^i\xi^j ] tQ_+^{2m}\dif x\dif t}_{IV} \\
&\underbrace{-\int_{0}^{t_0}\int_{\M}  [2Q_{ij}(\Delta \xi)^i\xi^j+4 Q_{ij;l} g^{kl}\xi^i_{;k}\xi^j+2Q_{ij}\xi^i_{;k}\xi^j_{;l}g^{kl}] tQ_+^{2m}\dif x\dif t}_V \\
:=&\uppercase\expandafter{\romannumeral1}+\uppercase\expandafter{\romannumeral2}+\uppercase\expandafter{\romannumeral3}+\uppercase\expandafter{\romannumeral4}+\uppercase\expandafter{\romannumeral5}.
\end{aligned}\label{24}
	\end{equation}
 The integrand vanishes outside $\{Q\ge 0\}\cap B(x_0,i_0)$, so we still can use the harmonic local coordinates.

We will bound 	$\uppercase\expandafter{\romannumeral1},\uppercase\expandafter{\romannumeral2},\uppercase\expandafter{\romannumeral4}$, $\uppercase\expandafter{\romannumeral5}$ in the form of $\uppercase\expandafter{\romannumeral3}$ and the left hand side respectively.
Using integration by parts,  and the fact that $Q_+=0$ at $t=0$,
 \begin{equation}\begin{aligned}
 	\uppercase\expandafter{\romannumeral1}=& -2m\int_{0}^{t_0}\int_{\M} tQ_+^{2m-1}|\nabla Q_+|^2\dif x\dif t\\
&-\frac{1}{2m+1}\int_MQ_+^{2m+1}(x,t_0)t_0\phi^2\dif x+\frac{1}{2m+1} \int_{0}^{t_0}\int_{\M}  Q_+^{2m+1}\dif x\dif t\\
\le &\frac{1}{2m+1} \int_{0}^{t_0}\int_{\M}  Q_+^{2m+1}\dif x\dif t.
 \end{aligned}\label{25}
 	 \end{equation}

 To bound other terms, we need the following estimate.
We view $H(\bullet,\xi)$ as a $(0,1)$-tensor, $|H(\bullet,\xi)|$ is its norm with respect to $g$.
To bound $|H(\bullet,\xi)|$, we choose an orthonormal basis $\{e_i\}$ that diagonalizes $H$  at any point $y\in B(x_0,i_0)$,
with $\frac{\pd}{\pd x^i}=P^k_i e_k$,  $H_{ij}=\lambda_i\delta_{ij}$, $\xi=\zeta^ie_i$.
\begin{equation}\al
Q_{ij}\xi^j &= H\left(\frac{\pd }{\pd x^i},\xi\right)=\sum\limits_{k=1}^n P_i^k\lambda_k \zeta^k \le
 \left(\sum\limits_{k=1}^n (\lambda_k \zeta^k)^2\right) ^{1/2}\left(\sum\limits_{k=1}^n (P_i^k)^2\right)^{1/2}
\\&\le \sqrt{C_1\xi^iQ_{ij}^2\xi^j }.
\eal
\label{26}
\end{equation}
where we used \eqref{07} and H\"older inequality. Hence
\begin{equation}
|H(\bullet,\xi)|=\left[g^{ij}H\left(\frac{\pd }{\pd x^i},\xi\right)H\left(\frac{\pd }{\pd x^j},\xi\right)\right]^{1/2}\le C_1  \sqrt{\xi^iQ_{ij}^2\xi^j}.
\label{27}
\end{equation}

For the second   term $II$  we still carry out integration by parts, use \eqref{13}, \eqref{21} and \eqref{27}, together with a recent result on  the sharp Li-Yau estimate for $- t \Delta \ln u$ in \cite{Z21}.

  \begin{eqnarray}\al
\uppercase\expandafter{\romannumeral2}=&\int_{0}^{t_0}\int_{\M}\left[ 2\langle \nabla\ln u,\nabla Q\rangle -4\nabla^k\ln u \xi^i_{;k} Q_{ij}\xi^j \right]tQ_+^{2m}  \dif x\dif t\\
 	&\le \int_{0}^{t_0} \int_{\M}\frac{2}{2m+1} (-t \Delta \ln u) Q_+^{2m+1}+  4t|\nabla\ln u||\nabla \xi| |H(\cdot,\xi)| Q_+^{2m}\dif x\dif t\\
 &\le \int_{0}^{t_0}\int_{\M} \bigg[ \frac{2{C(K,D_0,n,T)}}{2m+1} Q_+^{2m+1}  \\
 &+4\sqrt{C_3D_0^2+C_3T}\left( C_1C_2i_0^{-1}+\sqrt{C_1}\|\nabla\beta\|_\infty \right)C_1  \sqrt{\xi^iQ_{ij}^2\xi^j}Q_+^{2m} \bigg]  \dif x\dif t,
 \eal
 	 \label{28}\\
 	  \al
 	 \mbox{where}\qquad C(K,D_0,n,T)=\frac{n}{2}&+\sqrt{2n(n-1)K(1+(n-1)KT)(1+T)}D_0
 	 \\&+\sqrt{(n-1)K(1+(n-1)KT)(C_1+C_2K)T}.
 	 \eal
 	 \label{29}
 	\end{eqnarray}
 	
  Next,  using integration by parts on inner product  of tensors, we compute
    \be
    \al
    \uppercase\expandafter{\romannumeral4}= &\int_{0}^{t_0}\int_{\M} \Big[-2\xi^i R_{ik}Q_{kj}\xi^j-t R_{ij;l}(\nabla^l \ln u) \xi^i\xi^j +2 t R_{jl;i} (\nabla^l \ln u)\xi^i\xi^j \Big] tQ_+^{2m}\dif x\dif t \\
 	=&\int_{0}^{t_0}\int_{\M} -2\xi^i R_{ik}Q_{kj}\xi^j tQ_+^{2m}+ t^2R_{ij}(\xi^i\xi^j (\nabla^l\ln u)Q_+^{2m})_{;l}-2t^2 R_{jl}(\xi^i\xi^j (\nabla^l\ln u) Q_+^{2m})_{;i}\dif x\dif t\\
 	=& \int_{0}^{t_0}\int_{\M} -2\xi^i R_{ik}Q_{kj}\xi^j tQ_+^{2m}+ t^2R_{ij}\xi^i\xi^j (\De \ln u )Q_+^{2m}+2t R_{jl}\xi^i\xi^j Q_{il} Q_+^{2m}\\
 	&+\left( 2mR_{ij}(t\nabla^l\ln u)\xi^i\xi^j t Q_+^{2m-1}(\nabla_lQ_+)+ R_{ij}(t\nabla^l\ln u)\nabla_l(\xi^i\xi^j) tQ_+^{2m} \right)\\
 	&-\left(4m R_{jl}(t\nabla^l\ln u)\xi^i\xi^j t Q_+^{2m-1}(\nabla_iQ_+) + 2R_{jl}(t\nabla^l\ln u)\nabla_i(\xi^i\xi^j )tQ_+^{2m} \right)\dif x\dif t.
 \eal
 \ee Therefore
 \be
 \al
 IV	\le&\int_{0}^{t_0}  \int_{\M}  \bigg[(n-1)KC_1(C_3D_0^2+C_3T) Q_+^{2m}\\
 	&+\left( 2mR_{ij}(t\nabla^l\ln u)\xi^i\xi^j t Q_+^{2m-1}(\nabla_lQ_+)+ R_{ij}(t\nabla^l\ln u)\nabla_l(\xi^i\xi^j) tQ_+^{2m} \right)\\
 	&-\left(4m R_{jl}(t\nabla^l\ln u)\xi^i\xi^j t Q_+^{2m-1}(\nabla_iQ_+) + 2R_{jl}(t\nabla^l\ln u)\nabla_i(\xi^i\xi^j )tQ_+^{2m} \right)\bigg]\dif x\dif t.\\
 	 \eal
 \label{30}
 	 \ee	

In the above, we have used \eqref{21}, \eqref{22} and the bound of Ricci curvature $R_{ij}\le (n-1)K$ again. We further estimate:
  \begin{equation}\begin{aligned}
 &|\uppercase\expandafter{\romannumeral4}|\le \int_{0}^{t_0}\int_{\M} (n-1)KC_1(C_3D_0^2+C_3T)  Q_+^{2m} \dif x\dif t\\
 &+ \int_{0}^{t_0}\int_{\M}(n-1)KC_1\sqrt{C_3D_0^2+C_3T} \int_{0}^{t_0}\int_{\M} 6mt|\nabla Q_+|Q_+^{2m-1} +6t\| |\nabla\xi| \|_\infty  Q_+^{2m}
 \dif x\dif t\\
 \le& (n-1)KC_1\big[C_3D_0^2+C_3T+6T\sqrt{C_3D_0^2+C_3T}(C_1C_2i_0^{-1}+\sqrt{C_1}\|\nabla\beta\|_\infty)  \big]\int_{0}^{t_0}\int_{\M} Q_+^{2m}\dif x\dif t\\ &+m\int_{0}^{t_0}\int_Mt|\nabla Q_+|^2 Q_+^{2m-1}\dif x\dif t\\&+9m(n-1)^2K^2C_1^2T(C_3D_0^2+C_3T)\int_{0}^{t_0}\int_MQ_+^{2m-1} \dif x\dif t.
 \end{aligned}\label{31}
 \end{equation}	
where we use Young's inequality in the last row. The second term of (\ref{31}) can be absorbed by  negative terms in (\ref{25}),  then thrown away.

 Now we proceed to estimate $V$.
We use $H(\xi,\nabla_{\bullet}\xi)$ to represent the $(0,1)$-tensor $\bullet\mapsto H(\xi,\nabla_{\bullet}\xi)$.
We want to control the term $Q_{ij}\xi^i \xi^j_{;l}\nabla^lQ_+=\langle H(\xi,\nabla_{\bullet}\xi),\nabla Q_+\rangle$  in $\uppercase\expandafter{\romannumeral5}
$, for which special care is needed, in order to avoid the most negative eigenvalue of $Q_{ij}$. At any fixed  point $y \in  B(x_0,i_0)$, we  span $\frac{\xi(y)}{|\xi(y)|}$ into an orthonormal basis in $T_y M$ to
 construct a normal coordinates near $y$.  In doing so, we have an orthonormal basis $\{e_i\}$ at $y$ with $\nabla_{e_i}e_j=0$ and $\xi=\xi_1 e_1$  at $y$. Then
  \[
 \nabla_{e_k}\xi=\nabla_{e_k}(\xi_1e_1)=e_k(\xi_1)e_1+\xi_1\Gamma_{k1}^j e_{j}=e_k(\xi_1)e_1
 \]at $y$, $Q(y)=H_{11}\xi_1^2=H_{11}\eta^2\beta_1^2$.
$$\langle H(\xi,\nabla_{\bullet}\xi),\nabla Q_+\rangle(y)=H_{ij}\xi_i(\nabla_{e_k}\xi)_j Q_{+  ;k}(y)=H_{11}\xi_1 e_k(\xi_1)e_k(Q_+) (y)$$   In addition, from \eqref{10},  we have

\begin{equation}\al
|\langle H(\xi,\nabla_{\bullet}\xi),\nabla Q_+\rangle(y)|\le & |H_{11}||\eta\beta_1|e_k(\eta)e_k(Q_+)||\beta_1|(y)+|H_{11}||\eta^2\beta_1|| e_k(\beta_1)e_k (Q_+)|(y)\\
\le &|H_{11}||\eta\beta_1^2|C_2i_0^{-1}|\eta|^{1-\delta}|\nabla Q_+|(y) +|H_{11}\eta^2\beta_1^2||\beta|^{-1} |\nabla \beta||\nabla Q_+|(y)
\\\le& Q^{1-\frac{\delta}{2} }(C_1Q_0)^{\frac{\delta}{2}}C_2i_0^{-1}|\nabla Q_+|(y)+Q\|\nabla \beta\|_\infty \sqrt{C_1}|\nabla Q_+|(y).
\eal
\label{32}
\end{equation}
where we have used \eqref{10},  \eqref{12} and \eqref{13}.
By the arbitrariness of $y$,
we apply \eqref{32} to $\uppercase\expandafter{\romannumeral5}$, after integrating by parts of inner product of tensors and using \eqref{27}, we get
\begin{equation}\al
\uppercase\expandafter{\romannumeral5}&=-\int_{0}^{t_0}\int_{\M}  [2Q_{ij}(\Delta \xi)^i\xi^j+4 Q_{ij;l} g^{kl}\xi^i_{;k}\xi^j+2Q_{ij}\xi^i_{;k}\xi^j_{;l}g^{kl}] tQ_+^{2m}\dif x\dif t
\\&=\int_{0}^{t_0}\int_{\M} \big[ 4(2m)Q_{ij}\xi^i_{;l}\xi^j\nabla^l Q_+tQ_+^{2m-1}+2Q_{ij}\xi^i_{;k}g^{kl}\xi^j_{;l} tQ_+^{2m}+2Q_{ij}(\De \xi)^i\xi^j tQ_+^{2m} \big] \dif x\dif t\\
&\le \int_{0}^{t_0}\int_{\M}
\bigg[ 8m (C_1Q_0)^{\frac{\delta}{2}} C_2i_0^{-1} |\nabla Q_+|tQ_+^{2m-\frac{\delta}{2} }+8m\|\nabla_\beta\|_\infty \sqrt{C_1} |\nabla Q_+|tQ_+^{2m}
\\&\quad + 2Q_0(C_1C_2i_0^{-1}+\sqrt{C_1}\|\nabla\beta\|_\infty)^2tQ_+^{2m} \\
&\quad +2(C_2i_0^{-2}\sqrt{C_1}+(n-1)K\sqrt{C_1}+2C_2i_0^{-1}\|\nabla\beta\|_\infty) C_1\sqrt{\xi^iQ_{ij}^2\xi^j} tQ_+^{2m} \bigg] \dif x\dif t.
\eal
\ee Thus
\be
\al
V&\le \int_{0}^{t_0}\int_{\M} \bigg[
mt|\nabla Q_+|^2 Q_+^{2m-1}+ 8m (C_1Q_0)^\delta C_2^2i_0^{-2} T Q_+^{2m-\delta+1}+8m\|\nabla\beta\|_\infty^2C_1 TQ_+^{2m+1}\\
&\quad + 2Q_0(C_1C_2i_0^{-1}+\sqrt{C_1}\|\nabla\beta\|_\infty)^2TQ_+^{2m}\\
&\quad +2(C_2i_0^{-2}\sqrt{C_1}+(n-1)K\sqrt{C_1}+2C_2i_0^{-1}\|\nabla\beta\|_\infty) C_1\sqrt{\xi^iQ_{ij}^2\xi^j} TQ_+^{2m} \bigg] \dif x\dif t.
\eal\label{33}
\end{equation}
 In the above, we again used Young's inequality so that the first term on the right hand side of \eqref{33} can be  absorbed by negative terms  in \eqref{25}, then  discarded.

Substituting (\ref{33}), (\ref{31}) (\ref{28}) and (\ref{25}) into (\ref{24}), we have
\begin{equation}
	\al
		&2\int_{0}^{t_0}\int_{\M} \xi^iQ_{ij}^2\xi^jQ_+^{2m}\\
		\le &\int_{0}^{t_0}\int_{\M}   \frac{1}{2m+1}  Q_+^{2m+1} +\frac{2{C(K,D_0,n,T)}}{2m+1} Q_+^{2m+1}  \\
 &+\left(\frac{1}{2} \xi^iQ_{ij}^2\xi^j +8(C_3D_0^2+C_3T)\left( C_1C_2i_0^{-1}+\sqrt{C_1}\|\nabla\beta\|_\infty \right)^2C_1^2\right)Q_+^{2m}  \\
		 &+Q_+^{2m+1}+2C_1C_3D_0^2K Q_+^{2m} +2C_1C_3KT Q_+^{2m} +2(2n-1)C_1K TQ_0 Q_+^{2m} \dif x \dif t  + Y,
\eal
\ee where
\be
\al
Y \equiv &\int^{t_0}_{0} \int_{\M} \bigg[ (n-1)KC_1\big[C_3D_0^2+C_3T+6T\sqrt{C_3D_0^2+C_3T}(C_1C_2i_0^{-1}+\sqrt{C_1}\|\nabla\beta\|_\infty)  \big] Q_+^{2m}\\
&+9m(n-1)^2K^2C_1^2T^2(C_3D_0^2+C_3T)Q_+^{2m-1} \\
 &+ { 8m (C_1Q_0)^\delta C_2^2i_0^{-2} T Q_+^{2m-\delta+1}+8m\|\nabla\beta\|_\infty^2C_1 TQ_+^{2m+1} } \\
& + 2Q_0(C_1C_2i_0^{-1}+\sqrt{C_1}\|\nabla\beta\|_\infty)^2TQ_+^{2m} \\
&+\left(\frac{1}{2}\xi^iQ_{ij}^2\xi^j+2(C_2i_0^{-2}\sqrt{C_1}+(n-1)K\sqrt{C_1}+2C_2i_0^{-1}\|\nabla\beta\|_\infty)^2 C_1 ^2T ^2\right)Q_+^{2m} \bigg] \dif x\dif t.
	\eal\label{34}
	\end{equation}

Define $F_m:=\left(\int_{0}^{t_0}\int_{\M} Q_+^{m} dxdt\right)^{1/m}$,  inequality \eqref{34} becomes
\begin{equation}
	\begin{aligned}
&F^{2m+2}_{2m+2}=\int_{0}^{t_0}\int_{\M} Q_+^{2m+2} \dif x\dif t \le\quad C_1\int_{0}^{t_0}\int_{\M} \xi^iQ_{ij}^2\xi^jQ_+^{2m} \dif x\dif t
	\\ \le \quad&C_1\left[ \frac{1+2C(K,D_0,n,t)}{2m+1} +1+8m\|\nabla\beta\|_\infty^2C_1 T\right]F_{2m+1}^{2m+1}
	\\+&{ C_18m (C_1Q_0)^\delta C_2^2i_0^{-2} T F^{2m-\delta+1}_{2m-\delta+1} }\\
+& C_1\Big[ 8(C_3D_0^2+C_3T)(C_1C_2i_0^{-1}+\sqrt{C_1}\|\nabla\beta\|_\infty)^2C_1^2+2C_1C_3D_0^2K +2C_1C_3KT \\
&+2(2n-1)C_1KTQ_0\\&+(n-1)KC_1\big[C_3D_0^2+C_3T+6T\sqrt{C_3D_0^2+C_3T}(C_1C_2i_0^{-1}+\sqrt{C_1}\|\nabla\beta\|_\infty)  \big]\\
&+2Q_0(C_1C_2i_0^{-1}+\sqrt{C_1}\|\nabla\beta\|_\infty)^2T\\
&+2(C_2i_0^{-2}\sqrt{C_1}+(n-1)K\sqrt{C_1}+2C_2i_0^{-1}\|\nabla\beta\|_\infty)^2 C_1 ^2T ^2  \Big] F_{2m}^{2m}\\
 	 +&C_1 9m(n-1)^2K^2C_1^2T^2(C_3D_0^2+C_3T)F_{2m-1}^{2m-1}.
	\end{aligned}\label{35}
\end{equation}

Multiplying both sides with $\left(\int_{0}^{t_0}\int_{\M}1\dif x\dif t\right)^{-1}$ and using H\"older inequality\begin{equation}
	F_{2m+k}^{2m+k}\le F_{2m+2}^{2m+k}\left(\int_{0}^{t_0}\int_{\M}1\dif x\dif t\right)^{\frac{2-k}{2m+2}},\quad k=-1,0,1,
	\label{36}
	\end{equation}
then letting $X_m:=F_{2m+2}\left(\int_{0}^{t_0}\int_{\M}1\dif x\dif t\right)^{-\frac{1}{2m+2}}$,
we deduce,  using
\[
Q_0^{\delta} Q_+^{2m-\delta+1} \le Q_0 Q_+^{2m}
\]that
\begin{equation}
\al
		&X_m^{2m+2}\\
		\le& C_1\left[ \frac{1+2C(K,D_0,n,T)}{2m+1} +1+8m\|\nabla\beta\|_\infty^2C_1 T+8mC_1Q_0 C_2^2i_0^{-2} T \right]X_{m}^{2m+1}\\
+& C_1\Big[ 8(C_3D_0^2+C_3T)(C_1C_2i_0^{-1}+\sqrt{C_1}\|\nabla\beta\|_\infty)^2C_1^2+2C_1C_3D_0^2K +2C_1C_3KT \\
&+2(2n-1)C_1KTQ_0\\&+(n-1)KC_1\big[C_3D_0^2+C_3T+6T\sqrt{C_3D_0^2+C_3T}(C_1C_2i_0^{-1}+\sqrt{C_1}\|\nabla\beta\|_\infty)  \big]\\
&+2Q_0(C_1C_2i_0^{-1}+\sqrt{C_1}\|\nabla\beta\|_\infty)^2T\\
&+2(C_2i_0^{-2}\sqrt{C_1}+(n-1)K\sqrt{C_1}+2C_2i_0^{-1}\|\nabla\beta\|_\infty)^2 C_1 ^2T ^2  \Big] X_{m}^{2m}\\
 	 +&C_1 9m(n-1)^2K^2C_1^2T^2(C_3D_0^2+C_3T)X_m^{2m-1}.
	\eal
	\label{37}
\end{equation}

On the last term of right hand side of \eqref{37} , we have  coefficients growing as fast as $m$, hence we cannot directly iterate the inequality. Instead we shall obtain a $L^m$ bound for a fixed large number $m$.

 We can divide $X_m^{2m}$ on both sides since $X_m>0$. See that $X_m^2-\alpha X_m-\beta =\gamma X_m^{-1}$ has only one positive solution, so $X_m^3-\alpha X_m^2-\beta X_m-\gamma=0$ has only one positive root, for $\alpha, \beta,\gamma>0$. With the help of formula of roots for cubic equation, \eqref{37} implies that
\begin{equation}
	X_m\le C_1(m,K,n,D_0,i_0)+\sqrt{Q_0}C_2(m,K,n,D_0,i_0)+\sqrt[3]{Q_0}C_3(m,K,n,D_0,i_0).
\label{38}
\end{equation}
when $m$ is large enough and $T$ is fixed. The main point is that the power on $Q_0$ on the right hand side is less than $1$.

\medskip

{\it Step 3. $L^\infty$ bound.}

 To this end we will use a modified version of Moser's iteration, taking advantage of the good nonlinear term. Comparing with the standard Moser's iteration, one extra hurdle that needs to be overcome is the appearance of the term $Q_0$ on the right hand side. Recall that our goal is to find an upper bound for $Q_0$. So we should expect it to be on the left side after iteration. Fortunately, in the end the power on $Q_0$ on the right hand side is less than $1$.

Putting $Z=(Q-\alpha)_+$, with $\alpha$ being a constant to be given in \eqref{defalpha} below. For any positive integer $m$, we employ  $Z^{2m+1}$ as a test function on (\ref{23}).
By standard  integration argument, similar to the beginning of Step 2, we deduce that
 \begin{equation}\begin{aligned}
		&\int_{\M} \left[ \frac{2m+1}{(m+1)^2} |\nabla Z^{m+1}|^2+Z^{2m+1}Z_t \right] \dif x\\
		&+\int_{\{Q-\alpha\ge 0\}} \frac{2\xi^iQ^{ 2} _{ij}\xi^j-Q- 2C_1C_3D_0^2K }{t}Z^{2m+1} \dif x
\\
		&\le  \int_{\M} \bigg[ (2C_1C_3K + 2(2n-1)C_1KQ_0) Z^{2m+1} + \underbrace{2(\nabla_k Q-2Q_{ij}\xi^i \xi^j_{;k})(\nabla^k\ln u) Z^{2m+1}}_{T_1}
\\&\qquad\underbrace{-2\xi^iR_{ik}Q_{kj}\xi^jZ^{2m+1}-t( R_{ij;l}- R_{jl;i}-R_{il;j})(\nabla^l \ln u)\xi^i\xi^j Z^{2m+1}}_{T_2}
\\&\qquad \underbrace{-\left(2Q_{ij}(\Delta \xi)^i\xi^j+4 Q_{ij;l} g^{kl}\xi^i_{;k}\xi^j+2Q_{ij}\xi^i_{;k}\xi^j_{;l}g^{kl} \right) Z^{2m+1}}_{T_3} \bigg] \dif x.
	\end{aligned}\label{39}
	\end{equation}

Let us write \begin{equation}\al
	&T_1=2(\nabla_k Q-2Q_{ij}\xi^i \xi^j_{;k})(\nabla^k\ln u) Z^{2m+1},\\
	&T_2=-2\xi^iR_{ik}Q_{kj}\xi^jZ^{2m+1}-t( R_{ij;l}- R_{jl;i}-R_{il;j})(\nabla^l \ln u)\xi^i\xi^jZ^{2m+1},\\
	&T_3=-\left(2Q_{ij}(\Delta \xi)^i\xi^j+4 Q_{ij;l} g^{kl}\xi^i_{;k}\xi^j+2Q_{ij}\xi^i_{;k}\xi^j_{;l}g^{kl}\right)Z^{2m+1}.
	\eal
	\label{40}
\end{equation}

 We aim to bound the integrals of $T_i$, $i=1, 2, 3,$ in a way so that we can carry out Moser's iteration.

For $T_1$,we integrate in $x$ and apply \eqref{13}, \eqref{21}, \eqref{27} and the sharp Li-Yau estimates \cite{Z21}:
\begin{equation}\al
	I_1  \equiv \int _{\M} T_1\dif x &\le \int _{\M} \frac{1}{m+1}  (-\De\ln u)Z^{2m+2} +4|H(\bullet,\xi)||\nabla\xi||t\nabla\ln u|\frac{Z^{2m+1}}{t}\dif x \\
	&\le \int _{\M}  \bigg[\frac{C(K,D_0,n,T)}{m+1}\frac{Z^{2m+2}}{t} \\
	&\qquad+4C_1\sqrt{\xi^iQ_{ij}^2\xi^j}(C_1C_2i_0^{-1}+\sqrt{C_1}\|\nabla\beta\|_\infty)\sqrt{C_3D_0^2+C_3T}\frac{Z^{2m+1}}{t} \bigg] \dif x.
	\eal\label{41}
\end{equation}

 For $T_2$, after integration by parts, we deduce by Young's inequality and \eqref{22} that:
 \begin{equation}
 	\al
 	I_2  \equiv &\int _{\M} T_2\dif x\\
 	=&\int _{\M} \bigg[ (R_{ij}\xi^i\xi^j)(t\De  \ln u)Z^{2m+1}+\frac{2m+1}{m+1}R_{ij}\xi^i\xi^j(t\nabla^l\ln u )(\nabla_l Z^{m+1})Z^{m}
 	\\&+2R_{ij}(t\nabla^l\ln u)(\nabla_l \xi^i) \xi^jZ^{2m+1}
 	-\frac{2(2m+1)}{m+1} R_{jl}\xi^i\xi^j(t\nabla ^l \ln u)(\nabla_i Z^{m+1})Z^{m}
 	\\&-2R_{jl}(t\nabla^l\ln u)(\nabla_i\xi^i)\xi^jZ^{2m+1}-2R_{jl}(t\nabla^l\ln u)\xi^i(\nabla_i\xi^j)Z^{2m+1} \bigg] \dif x
 	\\ \le& \int _{\M} \bigg((n-1)KC_1(C_3D_0^2+C_3T)\frac{Z^{2m+1}}{t}
 	\\&+6(n-1)KC_1\sqrt{C_3D_0^2+C_3T}(C_1C_2i_0^{-1}+\sqrt{C_1}\|\beta\|_\infty) Z^{2m+1}
 	\\& +\frac{m}{2(m+1)^2} \left|\nabla Z^{m+1}\right|^2+\left[ \frac{9}{2}(n-1)^2K^2C_1^2(C_3D_0^2+C_3T)\frac{(2m+1)^2}{m}\right] Z^{2m} \bigg) \dif x. 	\eal
 	\label{42}
 \end{equation}

 For $T_3$, we use again \eqref{32}, \eqref{27} and Young's inequality to obtain:
 \begin{equation}\al
 I_3\equiv &\int _{\M} T_3\dif x \\
 	\le& \int_{\M} \bigg( \frac{4(2m+1)}{m+1} \Big[(Z+\alpha)^{1-\frac{\delta}{2}}(C_1Q_0)^{\delta/2} C_2i_0^{-1}|\nabla Z^{m+1}|Z^{m}+(Z+\alpha)\|\nabla \beta \|_\infty \sqrt{C_1}|\nabla Z^{m+1}|Z^m \Big]
 	\\&+2Q_0 (C_1C_2i_0^{-1}+\sqrt{C_1}\|\nabla\beta\|_\infty)Z^{2m+1}\\&+2(C_2i_0^{-2}\sqrt{C_1}+(n-1)K\sqrt{C_1}+2C_2i_0^{-1}\|\nabla\beta\|_\infty) C_1\sqrt{\xi^iQ_{ij}^2\xi^j} Z^{2m+1} \bigg) \dif x
 	\\\le& \int_{\M} \bigg[ \frac{m}{2(m+1)^2} \left|\nabla Z^{m+1}\right|^2+ 8\frac{(2m+1)^2}{m}(Z+\alpha)^{2-\delta}(C_1Q_0)^\delta Z^{2m}+8\frac{(2m+1)^2}{m} (Z+\alpha)^2C_1 \|\nabla\beta\|^2Z^{2m}
 		\\&+2Q_0 (C_1C_2i_0^{-1}+\sqrt{C_1}\|\nabla\beta\|_\infty)Z^{2m+1}
 		\\&+\frac{\xi^iQ_{ij}^2\xi^j }{t}Z^{2m+1}+(C_2i_0^{-2}\sqrt{C_1}+(n-1)K\sqrt{C_1}+2C_2i_0^{-1}\|\nabla\beta\|_\infty)^2 C_1^2TZ^{2m+1} \bigg] \dif x
 	\eal\label{43}
 \end{equation}

Substitute \eqref{41}, \eqref{42} and \eqref{43} into \eqref{39}, we have
 \begin{equation}\begin{aligned}
&\frac{1}{m+1}\int_{\M} \big[|\nabla  Z^{m+1}|^2+\frac{1}{2}\left(Z^{2m+2}\right)_t \big] \dif x+\mathrm I_4\\
\le & \int_{\M} \left[ 8\frac{(2m+1)^2}{m}Q^{2-\delta}(C_1Q_0)^\delta Z^{2m}+8\frac{(2m+1)^2}{m} (Z+\alpha)^2C_1\|\nabla\beta\|^2Z^{2m} \right]  \dif x
\\&+ \int_{\M} \Big[ 2C_1C_3K + 2(2n-1)C_1KQ_0 +6(n-1)KC_1\sqrt{C_3D_0^2+C_3T}(C_1C_2i_0^{-1}+\sqrt{C_1}\|\beta\|_\infty) \\
&+ 2Q_0(C_1C_2i_0^{-1}+\sqrt{C_1}\|\nabla\beta\|_\infty)
 + (\sqrt{C_1}C_2i_0^{-2}+(n-1)K\sqrt{C_1}+2C_2i_0^{-1}\|\nabla\beta\|_\infty)^2C_1^2T
\Big]Z^{2m+1}  \dif x
\\
&+\int_{\M} \left[ \frac{9}{2}(n-1)^2K^2C_1^2(C_3D_0^2+C_3T)\frac{(2m+1)^2}{m}+\frac{16(2m+1)^2}{m} (C_1C_2i_0^{-1}+\|\nabla\beta\|_\infty)^2 C_1nQ_0\alpha\right] Z^{2m}\dif x.
\end{aligned}\label{44}
\end{equation}
Here \begin{equation}\al
	\rm I_4=&\int_{\{Z\ge 0\}}\frac{Z^{2m+1}}{t} \Big[ 2\xi^iQ^{ 2} _{ij}\xi^j-(Z+\alpha)- 2C_1C_3D_0^2K  -\frac{C(K,D_0,n,t)}{m+1} Z
	\\&- 4C_1\sqrt{\xi^iQ_{ij}^2\xi^j}(C_1C_2i_0^{-1}+\sqrt{C_1}\|\nabla\beta\|_\infty)\sqrt{C_3D_0^2+C_3T}
	\\&-(n-1)KC_1(C_3D_0^2+C_3T) -\xi^iQ_{ij}\xi^j
	 \Big] \dif x
	 \\ =& \int_{\{Z\ge 0\}}\frac{Z^{2m+1}}{t} \Big[ \frac{\xi^iQ^{ 2}_{ij}\xi^j }{2}+ \frac{\xi^iQ^{ 2}_{ij}\xi^j }{4}-\left(1+\frac{C(K,D_0,n,T)}{m+1}\right) Z
	\\&+\frac{\xi^iQ^2_{ij}\xi^j }{4}- 4C_1(C_1C_2i_0^{-1}+\sqrt{C_1}\|\nabla\beta\|_\infty)\sqrt{C_3D_0^2+C_3T}  \sqrt{\xi^iQ_{ij}^2\xi^j}
	\\&-\alpha- 2C_1C_3D_0^2K-(n-1)KC_1(C_3D_0^2+C_3T) 	 \Big] \dif x
	\\ \ge &\int_{\{Z\ge 0\}} \frac{Z^{2m+1}}{t} \Big[ \frac{1}{2} C_1^{-1} \alpha ^2-C_1\left(1+\frac{C(K,D_0,n,T)}{m+1}\right)^2
	\\&- 16C_1^2(C_1C_2i_0^{-1}+\sqrt{C_1}\|\nabla\beta\|_\infty)^2(C_3D_0^2+C_3T)
	\\&-\alpha- 2C_1C_3D_0^2K-(n-1)KC_1(C_3D_0^2+C_3T)
	\Big] \dif x\\
	\ge & 0,
	\eal\label{45}
\end{equation}
if we choose
\be
\lab{defalpha}
\al\alpha=C_1+C_1\Bigg[1+2C_1\Big[C_1\left(1+\frac{C(K,D_0,n,t)}{m+1}\right)^2
	\\+16C_1^2(C_1C_2i_0^{-1}+\sqrt{C_1}\|\nabla\beta\|_\infty)^2(C_3D_0^2+C_3T)
	\\+ 2C_1C_3D_0^2K+(n-1)KC_1(C_3D_0^2+C_3T) \Big]\Bigg]^{1/2}.\eal
\ee  In the above we have used the inequality that
\[
\xi^iQ^{ 2}_{ij}\xi^j  |\xi|^2 \ge \alpha^2, \quad \text{when} \quad Q =\xi^i Q_{ij}\xi^j \ge \alpha,
\]which can be verified by diagonalizing $(Q_{ij})$ point-wise and using Cauchy-Schwarz inequality.

Using the lower bounds on curvature and volume of unit balls, it is well known that we have the Sobolev inequality: for $\mu=n/(n-2)$,
\begin{equation}
\int_{\M}\left|\nabla (Z^{m+1})\right|^2\dif x\ge  C_{S} v_0^{2/n}  \left(\int_{\M}\left|Z^{m+1}\right|^{2\mu} \dif x \right)^{\frac{1}{\mu}}
\label{46}
\end{equation}
with constant $C_S$ depending on $n$ and $K$.  See Chapter 18 of \cite{Lib} e.g.

Substitute (\ref{45}), (\ref{46}) into (\ref{44}) and integrate on $t \in [0, s] \subset [0, t_0]$, we deduce, since $Z=0$ when $t=0$, that
 \begin{equation}\begin{aligned}
 	& C_{S} v_0^{2/n}  \int_0^{s} \left(\int_{\M}\left|Z^{m+1}\right|^{2\mu}\dif x\right)^{\frac{1}{\mu}}\dif t+\frac{1}{2}\int_{\M} Z(x,s)^{2m+2}\dif x\\
 		\le &  \int_0^{s} \left(\int_{\M}\left|\nabla (Z^{m+1})\right|^2\dif x\right)^{\frac{1}{\mu}}\dif t+\frac{1}{2}\int_0^{s} \frac{\partial}{\partial t}\left(\int_{\M} Z^{2m+2}\dif x\right)\dif t\\
 		\le  &(m+1)\int_0^s\int_{\M} 8\frac{(2m+1)^2}{m}C_1Q_0 (Z^{2m+1}+\alpha Z^{2m})+8\frac{(2m+1)^2}{m} C_1\|\nabla\beta\|^2(Z^{2m+2}+2\alpha Z^{2m+1}+\alpha^2 Z^{2m})
\\&+ \Big[ 2C_1C_3K + 2(2n-1)C_1KQ_0 +6(n-1)KC_1\sqrt{C_3D_0^2+C_3T}(C_1C_2i_0^{-1}+\sqrt{C_1}\|\beta\|_\infty) \\
&+ 2Q_0(C_1C_2i_0^{-1}+\sqrt{C_1}\|\nabla\beta\|_\infty)
 + (\sqrt{C_1}C_2i_0^{-2}+(n-1)K\sqrt{C_1}+2C_2i_0^{-1}\|\nabla\beta\|_\infty)^2C_1^2T
\Big]Z^{2m+1}
\\
&+\left[ \frac{9}{2}(n-1)^2K^2C_1^2(C_3D_0^2+C_3T)\frac{(2m+1)^2}{m}+\frac{16(2m+1)^2}{m} (C_1C_2i_0^{-1}+\|\nabla\beta\|_\infty)^2 C_1nQ_0\alpha\right] Z^{2m}\dif x\dif t.
	\end{aligned}\label{47}
 		 	\end{equation}
 To reach this inequality, we used $Q_0^{\delta} Q_+^{2m-\delta+1} \le Q_0 Q_+^{2m}$ again.

On one hand, we have \begin{equation}
	\mbox{Right hand side of (\ref{47})} \ge \frac{1}{2}\sup\limits_{0\le s\le t_0}\int_{\M} Z(x,s)^{2m+2}\dif x:=\frac{1}{2}A.
 \label{48}
 \end{equation}

On the other hand, setting $s=t_0$ in (\ref{47}) gives us
\begin{equation}
\al
	(\mbox{Right hand side of (\ref{47}) with } s=t_0) &\ge  C_{S} v_0^{2/n} \int_{0}^{t_0} \left(\int_{\M} \left| Z^{m+1}\right|^{2\mu}\dif x\right)^{\frac{1}{\mu}}\dif t \\
&	:= C_{S} v_0^{2/n}  B,
\eal\label{49}
\end{equation}
By H\"older inequality, we have
\begin{equation}
 	\int_{0}^{t_0} \int_{\M}\left|Z^{m+1}\right|^{\frac{2(2\mu-1)}{\mu}}\dif x\dif t\le A^{\frac{\mu-1}{\mu}} B.\label{50}
 \end{equation}

 Multiplying (\ref{48}) to the $\frac{\mu-1}{2\mu-1}$th power  (\ref{49}) to the $\frac{\mu}{2\mu-1}$th power and using (\ref{50}), we  infer
\begin{equation}
	(\mbox{Right hand side of (\ref{47}) with } s=t_0)\ge \left(\frac{1}{2}\right)^{\frac{\mu-1}{2\mu-1}}  \left( C_{S} v_0^{2/n} \right)^{\frac{\mu}{2\mu-1}} \left( \int_{0}^{t_0} \int_{\M}\left|Z^{m+1}\right|^{\frac{2(2\mu-1)}{\mu}}\dif x\dif t\right)^{\frac{\mu}{2\mu-1}}.
	\label{51}
\end{equation}

We reorganize (\ref{51}) using Cauchy-Schwarz inequality \begin{equation}
	\begin{aligned}
		&\left( \int_{0}^{t_0} \int_{ \M}\left|Z\right|^{\frac{(2m+2)(2\mu-1)}{\mu}}\dif x\dif t\right)^{\frac{\mu}{2\mu-1}}\\
		\le &  \left(2\right)^{\frac{\mu-1}{2\mu-1}}  \left( C_{S} v_0^{2/n} \right)^{-\frac{\mu}{2\mu-1}} (m+1) \int_0^{t_0}\int_{\M} 8\frac{(2m+1)^2}{m} C_1\|\nabla\beta\|^2 Z^{2m+2} \dif x\dif t.
\\&+ \Big[ 2C_1C_3K + 2(2n-1)C_1KQ_0 +6(n-1)KC_1\sqrt{C_3D_0^2+C_3T}(C_1C_2i_0^{-1}+\sqrt{C_1}\|\beta\|_\infty) \\
&+ 2Q_0(C_1C_2i_0^{-1}+\sqrt{C_1}\|\nabla\beta\|_\infty)
 + (\sqrt{C_1}C_2i_0^{-2}+(n-1)K\sqrt{C_1}+2C_2i_0^{-1}\|\nabla\beta\|_\infty)^2C_1^2T
 \\&+8\frac{(2m+1)^2}{m}C_1Q_0 + 16\alpha\frac{(2m+1)^2}{m} C_1\|\nabla\beta\|^2
\Big]
\\&\cdot\left(\int_0^{t_0}\int_{\M}Z^{2m+2}\dif x\dif t\right)^{\frac{2m+1}{2m+2}}\left(\int_0^{t_0}\int_{\M}1\dif x\dif t\right)^{\frac{1}{2m+2}}
\\
&+\Big[ \frac{9}{2}(n-1)^2K^2C_1^2(C_3D_0^2+C_3T)\frac{(2m+1)^2}{m}+\frac{16(2m+1)^2}{m} (C_1C_2i_0^{-1}+\|\nabla\beta\|_\infty)^2 C_1nQ_0\alpha
\\&+ 8\frac{(2m+1)^2}{m}C_1Q_0 \alpha +8\frac{(2m+1)^2}{m} C_1\|\nabla\beta\|^2\alpha^2\Big]
\\
&\cdot\left(\int_0^{t_0}\int_{\M}Z^{2m+2}\dif x\dif t\right)^{\frac{2m}{2m+2}}
 \left(\int_0^{t_0}\int_{\M}1\dif x\dif t\right)^{\frac{2}{2m+2}}.
	\end{aligned}\label{52}
\end{equation}

Let us set $m=m_k=\frac{m_0}{2}\left(\frac{2\mu-1}{\mu}\right)^k-1$
and define  the sequence of numbers
 $$A_{p}:=\left(\int_{0}^{t_0}\int_{\M}	 Z^p\dif x\dif t\right)^\frac{1}{p}$$
We rewrite (\ref{52}) as
 	\begin{equation}\begin{aligned}
 			&A_{2m_{k+1}+2}^{2m_{k}+2}\\
 			\le& \left(2\right)^{\frac{\mu-1}{2\mu-1}}  {\left(C_s v_0^{2/n}\right)^{-\frac{\mu}{2\mu-1}} (m_k+1) }\Bigg[ 8\frac{(2m_k+1)^2}{m_k} C_1\|\nabla\beta\|^2 A^{2m_k+2}_{2m_k+2}
 			\\&+\Big[ 2C_1C_3K + 2(2n-1)C_1KQ_0 +6(n-1)KC_1\sqrt{C_3D_0^2+C_3T}(C_1C_2i_0^{-1}+\sqrt{C_1}\|\beta\|_\infty) \\
&+ 2Q_0(C_1C_2i_0^{-1}+\sqrt{C_1}\|\nabla\beta\|_\infty)
 + (\sqrt{C_1}C_2i_0^{-2}+(n-1)K\sqrt{C_1}+2C_2i_0^{-1}\|\nabla\beta\|_\infty)^2C_1^2T
 \\&+8\frac{(2m_k+1)^2}{m_k}C_1Q_0 + 16\alpha\frac{(2m_k+1)^2}{m_k} C_1\|\nabla\beta\|^2
\Big]A_{2m_k+2}^{2m_k+1}
 \left(\int_0^{t_0}\int_{\M}1\dif x\dif t\right)^{\frac{1}{2m_k+2}}
\\&+\Big[ \frac{9}{2}(n-1)^2K^2C_1^2(C_3D_0^2+C_3T)\frac{(2m_k+1)^2}{m_k}+\frac{16(2m_k+1)^2}{m_k} (C_1C_2i_0^{-1}+\|\nabla\beta\|_\infty)^2 C_1nQ_0\alpha
\\&+ 8\frac{(2m_k+1)^2}{m_k}C_1Q_0 \alpha +8\frac{(2m_k+1)^2}{m_k} C_1\|\nabla\beta\|^2\alpha^2\Big] A_{2m_k+2}^{2m_k} \left(\int_0^{t_0}\int_{\M}1\dif x\dif t\right)^{\frac{2}{2m_k+2}}\Bigg].
	\end{aligned}
 	\end{equation}

Since $(a+b+c)^p\le a^p+b^p+c^p$ for $0<p<1$ and $a,b,c>0$, we deduce that
 	 \begin{equation}\begin{aligned}
 	&A_{2m_{k+1}+2}\\
 			\le& \Bigg[\left(2\right)^{\frac{\mu-1}{2\mu-1}}  { \left(C_s v_0^{2/n}\right)^{-\frac{\mu}{2\mu-1}}} (m_k+1)\Bigg] ^{\frac{1}{2m_k+2}}\Bigg\{ \Bigg[ 8\frac{(2m_k+1)^2}{m_k} C_1\|\nabla\beta\|^2 \Bigg]^{\frac{1}{2m_k+2}}A_{2m_k+2}
 			\\& +\Bigg[ 2C_1C_3K + 2(2n-1)C_1KQ_0 +6(n-1)KC_1\sqrt{C_3D_0^2+C_3T}(C_1C_2i_0^{-1}+\sqrt{C_1}\|\beta\|_\infty) \\
&+ 2Q_0(C_1C_2i_0^{-1}+\sqrt{C_1}\|\nabla\beta\|_\infty)
 + (\sqrt{C_1}C_2i_0^{-2}+(n-1)K\sqrt{C_1}+2C_2i_0^{-1}\|\nabla\beta\|_\infty)^2C_1^2T
 \\&+8\frac{(2m_k+1)^2}{m_k}C_1Q_0 + 16\alpha\frac{(2m_k+1)^2}{m_k} C_1\|\nabla\beta\|^2
\Bigg]^{\frac{1}{2m_k+2}}A_{2m_k+2}^{\frac{2m_k+1}{2m_k+2}}
 \left(\int_0^{t_0}\int_{\M}1\dif x\dif t\right)^{\frac{1}{(2m_k+2)^2}}
 \\&+\Bigg[ \frac{9}{2}(n-1)^2K^2C_1^2(C_3D_0^2+C_3T)\frac{(2m_k+1)^2}{m_k}+\frac{16(2m_k+1)^2}{m_k} (C_1C_2i_0^{-1}+\|\nabla\beta\|_\infty)^2 C_1nQ_0\alpha
\\&+ 8\frac{(2m_k+1)^2}{m_k}C_1Q_0 \alpha +8\frac{(2m_k+1)^2}{m_k} C_1\|\nabla\beta\|^2\alpha^2\Bigg]^{\frac{1}{2m_k+2}} A_{2m_k+2}^{\frac{2m_k}{2m_k+2}} \left(\int_0^{t_0}\int_{\M}1\dif x\dif t\right)^{\frac{2}{(2m_k+2)^2}}\Bigg\}
 	\end{aligned}
 	\end{equation}

 If $A_{2m_k+2}\le 1$, then $A_{2m_{k+1}+2}\le  (Cm_k^2+Cm_k^2Q_0)^{\frac{1}{2m_k+2}}$. If $A_{2m_k+2}\ge  1$, then $A_{2m_k+2}^{\frac{2m_k}{2m_k+2}} \le A_{2m_k+2}$. The right hand side is of the form $(Cm_k^2+Cm_k^2Q_0)^{\frac{1}{2m_k+2}}(1+(t_0\vol(\M))^{\frac{2}{(2m_k+2)^2}})A_{2m_k+2}$, $C=C(K,D_0,n,i_0,T)$.


Iterate this inequality by $k$  till a fixed $m_0$, we conclude that	\begin{equation}A_{2m_{k+1}+2}\le\prod\limits_{i=0}^k \Bigg[Cm_i^2Q_0+CQ_0+C\Bigg]^{\frac{1}{2m_i+2}} A_{2m_0+2}.
		\label{55}
	\end{equation}

 Since $\prod\limits_{i=0}^\infty \left( m_0\left(\frac{2\mu-1}{\mu}\right)^{i}\right)^{\frac{2i}{m_0\left(\frac{2\mu-1}{\mu}\right)^{i}}}<\infty $ and $(CQ_0)^{\sum\limits_{i=0}^\infty\frac{2}{m_0^2\left(\frac{2\mu-1}{\mu}\right)^{2i}}}\le C Q_0^{\frac{1}{4}}$ for $m_0$ large enough. We let $k$ go to infinity in (\ref{55}) and recall that for any fixed large $m_0$, we have proved in Step 2:
 \[
 A_{2m_0+2}\le C(m_0,K,D_0,n,i_0)(\sqrt{Q_0}+\sqrt[3]{Q_0}+1)(t_0\vol(\M))^{\frac{1}{2m_0+2}}.
 \] Thus we arrive at, since $t_0 \le T$,

 \begin{equation}
Q_0-\alpha\le\sup\limits_{\M\times [0,t_0]} Z\le  C(K,D_0,n,i_0, T) \sqrt{Q_0}\cdot (1+C(n,K) Q_0^{\frac{1}{4}})\le C(K,D_0,n,i_0, T)Q_0^{3/4}.
  \end{equation}
Recall $u=G(x, t+\delta, y)$, therefore we  conclude, after letting $\delta \to 0$ and renaming the main parameter, that

$$ -t\nabla_i\nabla_j\ln G(x, t, y)\le \gamma(T, n,K,v_0)g_{ij}(x).$$

For any other positive solution $u=u(x, t)$,  we can write $u=\int_{\M} G(x, t; y) u_0(y) dy$, where $u_0(y):=\lim\limits_{t\to 0}u(y, t)$ and follow the steps in \cite{LZ} p63-64 to prove the claimed bound. This completes the proof of Theorem \ref{main}.

\end{proof}

\clearpage	
\section{appendix}
Now we recall Hamilton's  estimates for  $\partial_t \ln u$ and $|\nabla \ln u|$  in \cite{Ha93},  which were used in the previous section. The following is a combination of Theorem 1.1 and Lemma 4.1 there.

\begin{theorem}(Hamilton)
Let $u$ be a positive solution with $u\le A$, and let $-(n-1)K$ be a non-positive lower bound on the Ricci curvatures of $\M$. Then\begin{equation*}
	t|\nabla u|^2\le (1+2(n-1)Kt)u^2\ln\frac{A}{u}.
\end{equation*}
\begin{equation*}
	t \partial_t \ln u\le 2(1+(n-1)Kt)\ln \frac{A}{u}.
\end{equation*}
\label{Ham}
\end{theorem}

In our case here, we apply it to $G(x,t+\delta;y)$  to deduce
\begin{equation*}
	(t+\delta)|\nabla\ln u|^2\le (1+2(n-1)K(t+\delta))\ln\frac{A}{u}.
\end{equation*}

In \cite{Z21} , it was proved, under conditions which are satisfied here,  that  if $u=G(x, t; y)$ is the heat kernel, then
\begin{equation*}
	\ln\frac{1}{t^{n/2} u}\le C(n)(1+K+Kt)+\frac{D_0^2}{2t}
\end{equation*}
and \begin{equation}
\lab{sharply}
\al
&	-t\Delta\ln u\le \frac{n}{2}+\sqrt{2n(n-1)K(1+(n-1)KT)(1+T)}D_0
 	 \\
 &\qquad +\sqrt{(n-1)K(1+(n-1)KT)(C_1+C_2K)T}.
 \eal
\end{equation} Afterwards the above sharp Li-Yau estimate for all positive solutions follows from a standard trick.

We now outline the proof of Theorem \ref{Ham}, for the completeness.
\begin{proof}
	Compute that\begin{equation*}\al
		\partial_t\left(\frac{|\nabla u|^2}{u}\right)=&\Delta\left(\frac{|\nabla u|^2}{u}\right)-\frac{2}{u}\left|\nabla_i\nabla_j u-\frac{\nabla_i u\nabla_j u}{u}\right|^2 -2R_{ij}\frac{\nabla_i u\nabla_j u}{u}\\
		\le &\Delta\left(\frac{|\nabla u|^2}{u}\right)-\frac{2}{nu}\left|\De  u-\frac{|\nabla u|^2}{u}\right|^2 +2(n-1)K\frac{|\nabla u|^2}{u}.
	\eal
	\end{equation*}
	and \begin{equation*}
		\partial_t\left(u\ln\frac{A}{u}\right)=\Delta\left(u\ln\frac{A}{u}\right) +\frac{|\nabla u|^2}{u}.
	\end{equation*}
	
	Consider $h_1=\varphi_1\frac{|\nabla u|^2}{u}-u \ln\frac{A}{u}$, where\begin{equation*}
		\varphi_1:=\frac{t}{1+2(n-1)Kt}.
	\end{equation*}
	
	Then $h_1(0)\le 0$ and $(\Delta-\partial_t)h_1\ge 0$. By maximum principle we have $h_1\le 0$ for all time.
	
	Consider $h_2=\varphi_2\left[\De u+\frac{|\nabla u|^2}{u}\right]-u \left[n+4\ln\frac{A}{u}\right]$, where\begin{equation*}
		\varphi_2:=\frac{e^{(n-1)Kt}-1}{(n-1)Ke^{(n-1)Kt}}.
	\end{equation*}
	
	Then $h_2(0)\le  0$ and \[(\Delta-\partial_t)h_2=\frac{2\varphi_2}{nu}\left[\De u-\frac{|\nabla u|^2}{u}\right]^2-\varphi^\prime\left[\De u-\frac{|\nabla u|^2}{u}\right]+2\frac{|\nabla u|^2}{u}\ge 0. \]whenever $h_2\ge 0$. By maximum principle we have $h_2\le 0$ for all time.
Thus \begin{equation*}
	t|\nabla u|^2\le (1+2(n-1)Kt)u^2\ln\frac{A}{u},\qquad
	t \partial_t \ln u\le 2(1+(n-1)Kt)\ln \frac{A}{u}.
\end{equation*}

\end{proof}

For \eqref{sharply},	
we consider $Y=Y(x,t)=-\De \ln u.$ $Y^+$ is a subsolution of the inequality:$$(\De-\pd_t)Y^++2\nabla Y^+\nabla\ln u\ge\frac{2}{n}(Y^+)^2-2(n-1)K|\nabla\ln u|^2.$$
Choosing $((t-\e )^+)^{2j+2}(Y^+)^{2j}$ as a test function and integrating, after suitable iteration, one can obtain the sharp Li-Yau gradient bound.

\section*{acknowledgements}
	The authors wish to thank Professors Bennett Chow and Xiaolong Li for their  help on the subject. The support of Simons Foundation grant 710364 is gratefully acknowledged.


\end{document}